\documentclass[11pt,bezier]{article}
\usepackage{amsmath}
\usepackage{amsfonts,amsthm,amssymb}
\usepackage{amsfonts}
\usepackage{graphics}
\usepackage{amssymb}
\newtheorem{lem}{Lemma}[section]
\newtheorem{thm}[lem]{Theorem}

\newtheorem{conj}{Conjecture}

\theoremstyle{definition}

\begin{document}
\title{Connectivity keeping edges of trees in 3-connected or 3-edge-connected graphs\footnote{The research is supported by National Natural Science Foundation of China (12261086).}}
\author{Qing Yang, Yingzhi Tian\footnote{Corresponding author. E-mail: tianyzhxj@163.com (Y. Tian).} \\
{\small College of Mathematics and System Sciences, Xinjiang
University, Urumqi, Xinjiang 830046, PR China}}

\date{}
\maketitle

\noindent{\bf Abstract } Hasunuma [J. Graph Theory 102 (2023) 423-435] conjectured that for any tree $T$ of order $m$, every $k$-connected (or $k$-edge-connected) graph $G$ with minimum degree at least $k+m-1$ contains a tree $T'\cong T$ such that $G-E(T')$ is still $k$-connected (or $k$-edge connected). Hasunuma verified this conjecture for $k\leq 2$. In this paper, we confirm this conjecture for $k=3$.


\section{Introduction}

The graphs considered in this paper are simple, finite and undirected. For graph-theoretical terminologies and notation not defined here, we follow \cite{Bondy}.

The \emph{connectivity} of $G$, denoted by $\kappa(G)$, is the minimum size of a vertex set $U\subseteq V(G)$ such that $G-U$ is disconnected or has only one vertex. The graph $G$ is said to be \emph{$k$-connected} if $\kappa(G)\geq k$. The \emph{edge-connectivity} of $G$, denoted by $\kappa'(G)$, is the minimum size of an edge set $F$ such that $G-F$ is disconnected. The graph $G$ is said to be \emph{$k$-edge-connected} if $\kappa'(G)\geq k$.



In 1972, Chartrand, Kaigars and Lick \cite{Chartrand} showed that there is a redundant vertex in a $k$-connected graph by the following theorem.

\begin{thm} (Chartrand, Kaigars and Lick \cite{Chartrand}) Every $k$-connected graph $G$ with minimum degree at least $\lfloor\frac{3k}{2}\rfloor$ contains a vertex $v$ such that $\kappa (G-\{v\})\geq k$.
\end{thm}

In 2008, Fujita and Kawarabayashi \cite{Fujita} showed that there are two redundant adjacent vertices in a $k$-connected graph with minimum degree at least $\lfloor\frac{3k}{2}\rfloor+2$. Furthermore, they proposed a conjecture relating to redundant subgraph in a $k$-connected graph. In 2010, Mader \cite{Mader1} conjectured that there is a redundant tree in a $k$-connected graph.

\begin{conj} (Mader \cite{Mader1}) For any tree $T$ of order $m$, every $k$-connected graph $G$ with minimum degree at least $\lfloor\frac{3k}{2}\rfloor+m-1$ contains a tree $T'\cong T$ such that $\kappa (G-V(T'))\geq k$.
\end{conj}

In the same paper, Mader confirmed Conjecture 1 for $T$ to be a path. Mader \cite{Mader2} also proved that Conjecture 1 holds if $\delta(G)\geq2(k-1+m)^2+m-1$. Actually, the result in \cite{Diwan} implies the correctness of Conjecture 1 for $k=1$. When $k=2$, there are several partial results for $T$ to be some special trees (such as star, double-star, spider and caterpillar et al.), see  [6-7, 9, 12, 16-17] for references. In 2022, Hong and Liu \cite{Hong} confirmed Conjecture 1 for $k=2$ and $k=3$. For $k\geq4$, Conjecture 2 remains open. Concerning bipartite graphs, Luo, Tian and Wu showed that every bipartite graph $G$ with $\kappa (G)\geq k$ and $\delta(G)\geq k+m$ contains a path $P$ of order $m$ such that $\kappa (G-V(P))\geq k$. And they proposed the following conjecture.

\begin{conj} (Luo, Tian and Wu \cite{Luo}) For any tree $T$ with bipartition $X$ and $Y$, every $k$-connected bipartite graph $G$ with minimum degree at least $k+t$, where $t=max\{|X|,|Y|\}$, contains a tree $T'\cong T$ such that $\kappa (G-V(T'))\geq k$.
\end{conj}

In \cite{Zhang}, Zhang confirmed Conjecture 2 for $T$ to be caterpillar and $k\leq2$. Furthermore, we verified Conjecture 2 for $T$ to be a caterpillar and $k=3$, and for $T$ to be a spider and $k\leq 3$ in \cite{Yang1}. Recently, we also confirmed Conjecture 2 for $T$ to be an odd path for any $k$ in \cite{Yang2}.

Naturally, there are similar problems concerning deletion edges of a subgraph in the $k$-connected graphs or $k$-edge-connected graphs.

\begin{thm} (Halin \cite{Halin}) Every graph $G$ with $\kappa (G)\geq k$ and $\delta(G)\geq k+1$ contains an edge $e$ such that $\kappa (G-\{e\})\geq k$.
\end{thm}

\begin{thm} (Lick \cite{Lick}) Every graph $G$ with $\kappa'(G)\geq k$ and $\delta(G)\geq k+1$ contains an edge $e$ such that $\kappa'(G-\{e\})\geq k$.
\end{thm}

Analogous to Conjecture 1, Hasunuma recently proposed the following edge-version conjecture.



\begin{conj} (Hasunuma \cite{Hasunuma3}) For any tree $T$ of order $m$, every $k$-connected graph (respectively, $k$-edge-connected graph) $G$ with minimum degree at least $k+m-1$ contains a tree $T'\cong T$ such that $\kappa (G-E(T'))\geq k$ (respectively, $\kappa' (G-E(T'))\geq k$).
\end{conj}

In \cite{Hasunuma3}, Hasunuma verified Conjecture 3 for $k\leq2$. For $k\geq3$, Conjecture 3 remains open.

Motivated by the results above, we will study Conjecture 3 for $k=3$. In the next section,  we will introduce some definitions and lemmas, which will be used to prove the main results. In Section 3, we will give the main results, which confirm Conjecture 3 for $k=3$.

\section{Preliminaries}
Let $G$ be a graph with \emph{vertex set} $V(G)$ and \emph{edge set} $V(G)$. The $order$ of $G$ is $|V(G)|$. For a vertex $v\in V(G)$, the \emph{neighborhood} of $v$ in $G$, denoted by $N_{G}(v)$, is the set of vertices in $G$ adjacent to $v$. The degree $d_G(v)$ of $v$ in $G$ is $|N_{G}(v)|$. The \emph{minimum degree} $\delta(G)$  of $G$ is min$_{v\in V(G)}d_G(v)$. The \emph{maximum degree} $\Delta(G)$  of $G$ is max$_{v\in V(G)}d_G(v)$. For a subgraph $H\subseteq G$, define $\delta_{G}(H)=\min_{v\in V(H)}d_{G}(v)$. While $\delta(H)$ denotes the minimum degree of the graph $H$. The \emph{induced subgraph} of a vertex set $U$ in $G$, denoted by $G[U]$, is the graph with vertex set $U$, where two vertices in $U$ are adjacent if and only if they are adjacent in $G$. And $G-U$ is the induced graph $G[V(G)\backslash U)]$. The \emph{edge-induced subgraph} of an edge set $F$ in $G$, denoted by $G[F]$, is the subgraph of $G$ with edge set $F$ and vertex set
consisting of all ends of edges in $F$. And $G-F$ is the subgraph of $G$ with vertex set $V(G)$ and edge set $E(G)\backslash F$.

Let $T$ be a tree. The set of internal vertices, denoted by $V_{I}(T)$, is the vertices with degree at least two in $T$. The set of leaves, denoted by $V_{L}(T)$, is the vertices with degree one in $T$. The cardinalities of $V_{I}(T)$ and $V_{L}(T)$ are denoted by $I(T)$ and $L(T)$, respectively.

The following Lemmas will be used to construct trees in our proof.

\begin{lem}(Hasunuma and Ono \cite{Hasunuma1})
Let $T$ be a tree of order $m$  and let $S$ be a subtree of $T$. If a graph $G$ contains a subtree $S'\cong S $ such that $d_{G}(v)\geq m-1$ for any $v\in V(G)\setminus V(S')$ and for any $v\in V(S')$ with $d_{S}(\phi ^{-1}(v))<d_{T}(\phi ^{-1}(v))$, where $\phi$ is an isomorphism from $V(S)$ to $V(S')$, then $G$ contains a subtree $T'\cong T$ such that $S'\subseteq T'$.
 \end{lem}

\begin{lem}(Hasunuma and Ono \cite{Hasunuma1})
Let $T$ be a tree of order $m$ and let $S$ be a subtree obtained from $T$ by deleting some leaves of $T$. If a graph $G$ contains a subtree $S'\cong S $ such that $d_{G}(v)\geq m-1$ for any $v\in V(S')$ with $d_{S}(\phi ^{-1}(v))<d_{T}(\phi ^{-1}(v))$, where $\phi$ is an isomorphism from $V(S)$ to $V(S')$, then $G$ contains a subtree $T'\cong T$ such that $S'\subseteq T'$.
 \end{lem}

\section{Main results}
The $subdivision$ for an edge $uv$ in a graph $G$  is the deletion of $uv$ from $G$ and the addition of two edges $uw$ and $wv$ along with the new vertex $w$. The subdivision of the graph $G$ is derived from $G$ by a sequence of subdivisions for some edges in  $G$.

Let $G$ be a subdivision of some simple 3-connected graph. An \emph{ear} of $G$ is a maximal path $P$ whose each internal vertex has degree 2 in $G$. A \emph{$(u, v)$-path} $P$ is a path with ends $u$ and $v$. For any $u',v'\in V(P)$, $u'Pv'$ is the subpath of $P$ between $u'$ and $v'$. For two subsets $V_{1}$ and $V_{2}$ of $V (G)$, the $(V_{1}, V_{2})$-path is a path with one end in $V_{1}$, the other end in $V_{2}$, but internal vertices not in $V_{1}\cup V_{2}$. We use $(v, V_{2} )$-path for $(\{v\}, V_{2})$-path.
Denote by $[V_1,V_2]_G$ the set of edges connecting one vertex in $V_1$ and one vertex in $V_2$.

Theorems 1.2 and 1.3 imply that Conjecture 3 holds for $m=2$. So we assume $m\geq3$ in our results.

\begin{thm}
For any tree $T$ of order $m\ (\geq3)$, every 3-connected graph $G$ with $\delta(G)\geq m+2$ contains a subtree $T'\cong T$ such that $G-E(T')$ is still 3-connected.
\end{thm}

\noindent{\bf Proof.} By contradiction, assume the theorem is false. Since $\delta(G)\geq m+2$, there exists a tree $T'\cong T$ in $G$ by Lemma 2.1. Let $B$ be a subdivision of some simple 3-connected graph of $G-E(T')$. In fact, $B$ does exist, since $\delta(G-E(T'))\geq 3$. Furthermore, we choose $T'$ and $B$ such that $n(B)$ is as large as possible and $|V(B)|$ is as small as possible, where $n(B)=|\{v\in V(B)|d_{B}(v)\geq 3\}|$. And we denote $$V_{1}=\{v\in V(B)|d_{B}(v)\geq 3\}$$ and $V_{2}=V(B)\backslash V_{1}$.

Suppose $V(G)\backslash V(B)=\emptyset$. Then $\delta(B)\geq m+2-\Delta(T)\geq m+2-(m-1)=3$ and $B$ is 3-connected, a contradiction. Thus, $V(G)\backslash V(B)\neq\emptyset$. We will complete the proof by a series of claims in the following.

\noindent{\bf Claim 1.} For any tree $T^{*}(\cong T)$ in $G$, if $G-E(T^{*})$ contains $B$, then $|N_{G-E(T^{*})}(v)\cap V(B)|\leq 3$ for each $v\in V(G)\backslash V(B)$.

By contrary, assume that $w$ is such a vertex satisfying $|N_{G-E(T^{*})}(w)\cap V(B)|\geq 4$.

Suppose no ear of $B$ contains $N_{G-E(T^{*})}(w)\cap V(B)$. Let $B_1$ be the subgraph obtained from $B$ by adding the vertex $w$ and the edge set $[w, V(B)]_{G-E(T^{*})}$.
Since $|N_{G-E(T^{*})}(w)\cap V(B)|\geq 4$, there are three disjoint $(w,V_{1})$-paths in $B_1$.
Thus $B_1$ is a subdivision of some simple 3-connected graph of $G-E(T^*)$.
But $n(B_1)> n(B)$, which contradicts to the choices of $T'$ and $B$.

Suppose there exists an ear $Q$ containing $N_{G-E(T^{*})}(w)\cap V(B)$. Let $Q=u_{1}u_{2}\cdots u_{t}$, $p=\min\{i:1\leq i\leq t ~and~u_{i}\in N_{G-E(T^{*})}(w)\}$ and $q=\max\{i:1\leq i\leq t~and~u_{i}\in N_{G-E(T^{*})}(w)\}$. Then by $|N_{G-E(T^{*})}(w)\cap V(B)|\geq 4$, we have $t\geq 4$ and $q-p\geq 3$. Let $B_2$ be the subgraph obtained from $B$ by deleting the vertex set $\{u_{p+1},\cdots,u_{q-1}\}$, and adding the vertex $w$ and edges $wu_p$, $wu_q$. Then $B_2$ is still a subdivision of some simple 3-connected graph in $G-E(T^*)$. But $n(B_2)= n(B)$ and $|V(B_2)|<|V(B)|$, also a contradiction.

\noindent{\bf Claim 2.} $|N_{G}(v)\cap V(B)|\leq 3+\Delta(T)$ for each $v\in V(G)\backslash V(B)$.

By Claim 1, we have $|N_{G-E(T')}(v)\cap V(B)|\leq 3$ for every $v\in V(G)\backslash V(B)$. Thus $|N_{G}(v)\cap V(B)|\leq 3+\Delta(T)$ is obtained.

\noindent{\bf Claim 3.} $|N_{G}(v)\cap V(B)|\leq 3$ for each $v\in V(G)\backslash V(B)$.

By contrary, assume that $x$ is such a vertex satisfying $|N_{G}(x)\cap V(B)|\geq 4.$ Let us consider the following two cases.

\noindent{\bf Case 1.} $T$ is a star.

If $V(B)\cup \{x\}\neq V(G)$, then there exists a vertex $y\in V(G)\backslash (V(B)\cup \{x\})$. Let $$T_1=G[\{ yy_{1},yy_{2},\cdots,yy_{m-1}\}],$$ where $\{y_{1},y_{2},\cdots,y_{m-1}\}\subseteq N_{G}(y)$. Then $T_1$ is isomorphic to $T'$ and $G-E(T_1)$ contains $B$. But $|N_{G-E(T_1)}(x)\cap V(B)|\geq 4$, which contradicts to Claim 1.

If $V(B)\cup \{x\}=V(G)$, then $x$ is the center of $T'$. Otherwise, $|N_{G-E(T')}(x)\cap V(B)|\geq m+2-1\geq 4$, which contradicts to Claim 1. Thus $\delta (B)\geq m+2-1\geq 4$ and $B$ is a 3-connected graph. Moreover, we have $|N_{G-E(T')}(x)\cap V(B)|\geq m+2-(m-1)=3$ by $d_{G}(x)\geq m+2$. Therefore, $G-E(T')$ is 3-connected, which contradicts to the assumption.

\noindent{\bf Case 2.} $T$ is not a star.

Since $m=I(T)+L(T)\geq I(T)+\Delta(T)$, we have $\delta(G-V(B))\geq m+2-(3+\Delta(T))\geq I(T)-1$ by Claim 2. Thus, we greedily construct a subtree $S_1\cong T-V_{L}(T)$ in $G-V(B)$ with an isomorphism $\phi$ from $V(T)\backslash V_{L}(T)$ to $V(S_1)$. By $\delta(G)\geq m+2$ and Lemma 2.2, we can extend $S_1$ to a tree $T_2\cong T$ in $G-E(B)$. Since $$|N_{G}(x)\backslash V(S_1)|\geq m+2-(I(T)-1)\geq\Delta(T)+3,$$ even if $x\in V(S_1)$ and $d_{T-V_{L}(T)}(\phi ^{-1}(x))<d_{T}(\phi ^{-1}(x))$, by doing the extension process from $S_1$ to $T_2$, we can leaving four edges joining $x$ and vertices in $B$ unused. But $|N_{G-E(T_2)}(x)\cap V(B)|\geq 4$, which contradicts to Claim 1. Thus Claim 3 holds.

\noindent{\bf Claim 4.} $V_{2}\neq \emptyset$.

By contrary, assume that $V_{2}=\emptyset$. Then $B$ is 3-connected. Since $G$ is 3-connected, for any $u\notin V(B)$, there exist three disjoint $(u,V(B))$-paths $P_{1}$, $P_{2}$, $P_{3}$ in $G$. We assume $P_{1}$, $P_{2}$, $P_{3}$ are such that $\Sigma_{i=1}^{3}|V(P_{i})|$ is as small as possible. Let $B_3$ be the subgraph obtained from $B$ by adding both the vertices and edges of $P_{1}$, $P_{2}$ and $P_{3}$.
Then $n(B_3)> n(B)$ and $B_3$ is still the subdivision of some simple 3-connected graph.

Suppose $\delta(G-V(B_3))\geq I(T)-1$. Then there exists a subtree $S_2\cong T-V_{L}(T)$ in $G-V(B_3)$. We can extend $S_2$ to a tree $T_3\cong T$ in $G-E(B_3)$ by Lemma 2.2. But $n(B_3)> n(B)$, which contradicts to the choices of $T'$ and $B$.

Suppose $\delta(G-V(B_3))\leq I(T)-2$. Let $z$ be a vertex of $G-V(B_3)$ with minimum degree. Then $$|N_{G}(z)\cap V(B_3)|\geq m+2-(I(T)-2)\geq 6$$ by $I(T)\leq m-2$. For $i=1,2,3$, let $z_{i}\in N_{G}(z)\cap V(B_3)$ such that $z_{i}\notin \{z_{1},\cdots,z_{i-1}\}$ and the distance from $z_{i}$ to $B$ in $B_3-\{z_{1},\cdots,z_{i-1}\}$ is as small as possible. Then there exist three disjoint $(z,V(B))$-paths $P'_{1}$, $P'_{2}$, $P'_{3}$ in $G[V(B_3)\cup \{z\}]$ using edges $zz_{1},zz_{2},zz_{3}$ such that $\Sigma_{i=1}^{3}|V(P'_{i})|<\Sigma_{i=1}^{3}|V(P_{i})|$, a contradiction. Thus Claim 4 holds.

By $\kappa(G)\geq3$ and Claim 4, there exists a $(V_{2},V(B))$-path $P$ such that the two ends of $P$ lie in no ear of $B$. In addition, we choose such a path such that $|V(P)|$ is as small as possible.
Let $B_4$  be the subgraph obtained from $B$ by adding both the vertices and edges of the path $P$. Then $n(B_4)>n(B)$ and $B_4$ is still the subdivision of some simple 3-connected graph.

Suppose $\delta(G-V(B_4))\geq I(T)-1$. Then there exists a subtree $S_3\cong T-V_{L}(T)$ in $G-V(B_4)$. We can extend $S_3$ to a tree $T_4\cong T$ in $G-E(B_4)$ by Lemma 2.2. But $n(B_4)> n(B)$, a contradiction.

Suppose $\delta(G-V(B_4))\leq I(T)-2$. Let $s$ be a vertex of $G-V(B_4)$ with minimum degree. Then $|N_{G}(s)\cap V(B_4)|\geq m+2-(I(T)-2)\geq 6$. Assume $P$ is a $(u,u')$-path such that $u\in V_{2}$.
By Claim 3, we have $$|N_{G}(s)\cap V(P)|\geq |N_{G}(s)\cap V(B_4)|-|N_{G}(s)\cap (V(B)\setminus\{u,u'\})|\geq 3.$$
By the choice of $P$, we have $|N_{G}(s)\cap V(P)|=3$, $|N_{G}(s)\cap (V(B)\setminus\{u,u'\})|=3$ and $N_{G}(s)\cap \{u,u'\}=\emptyset$. Let $s'\in N_{G}(s)\cap V(P) $ closest to $u$ in $P$. And let $s''\in N_{G}(s)\cap (V(B)\setminus\{u,u'\})$ such that $s''$ is not an end of the ear containing $u$ in $B$. Denote $P'=uPs'ss''$. Then $P'$ is a $(V_{2},B)$-path such that the two ends of $P'$ lie in no ear of $B$. But $|V(P')|<|V(P)|$, which contradicts to the choice of $P$. $\Box$

By a similar argument, we can obtain the desired result in 3-edge-connected graphs.

\begin{thm}
For any tree $T$ of order $m\ (\geq3)$, every 3-edge connected graph $G$ with $\delta(G)\geq m+2$ contains a subtree $T'\cong T$ such that $\kappa'(G-E(T'))\geq 3$.
\end{thm}

\end{document}